\documentclass[a4paper,pdftex]{article}
\usepackage{amsmath,amssymb,latexsym,amsfonts}

\usepackage{color}

\usepackage{euscript}
\newcommand{\x}{\times}

\newcommand{\diag}{\mathop{\rm diag}\nolimits}

\def\cqd{\ensuremath{\Box}}

\def\qed{\ensuremath{\blacksquare}}

\newcommand{\R}{\ensuremath{\mathbb R}}
\newcommand{\C}{\ensuremath{\mathbb C}}

\newcommand{\mi}{\mathrm{i}}

\newcommand{\al}{\ensuremath{\alpha}}

\newcommand{\la}{\ensuremath{\lambda}}

\newcommand{\si}{\ensuremath{\sigma}}

\newcommand{\Om}{\ensuremath{\Omega}}

\newcommand{\Cmn}{\ensuremath{{\mathbb C}^{m \times n}}}
\newcommand{\Cnn}{\ensuremath{{\mathbb C}^{n \times n}}}

\newcommand{\Cmuno}{\C^{m\times 1}}
\newcommand{\Cnuno}{\C^{n\times 1}}

\newcommand{\Cmmasn}{\C^{(m+n)\x (m+n)}}

\newcommand{\n}[1]{\boldsymbol{#1}}

\newcommand{\s}[1]{\ensuremath{\Lambda(#1)}}

\newcommand{\secsim}[2]{\ensuremath{#1_1,
\dots,#1_{#2}  }}  

\newcommand{\secdec}[2]{\ensuremath{{#1}_1\ge \cdots \ge {#1}_{#2}}}

\newcommand{\secudectres}[3]{\ensuremath{#1_1(#2)\ge #1_2(#2) \ge \cdots \ge #1_#3(#2)}}

\def\demo{\noindent{\sc Proof.} }

\def\recoigu#1#2#3{\ensuremath{#1 = #2,\dots,#3}}

\newcommand{\recoper}[3]{\ensuremath{#1 \in \{#2,\dots,#3\} }}

\newcommand{\norma}[1]{\lVert #1 \rVert}

\newcommand{\secuencia}[3]{#1_{#2},\ldots , #1_{#3}}

\newcommand{\sectres}[3]{#1_1(#2),\ldots,#1_#3(#2)}

\newcommand{\secudecuatro}[4]{\ensuremath{#1}_{#2}\left(#3\right) \ge \cdots \ge {#1}_{#4}\left(#3\right)}

\newcommand{\SecuIgualCuatro}[4]{\ensuremath{#1}_{#2}\left(#3
\right) = \cdots = {#1}_{#4}\left(#3\right)}

\newcommand{\SecuSumaCuatro}[4]{\ensuremath{#1}_{#2}\left(#3
\right) + \cdots + {#1}_{#4}\left(#3\right)}

\newcommand{\parcial}[2]{\frac{\partial {#1}}{\partial {#2}}}

\renewcommand{\ge}{\geqslant}
\renewcommand{\le}{\leqslant}

\usepackage[pdftex]{graphicx}
\usepackage{upref,ifthen}

\usepackage[numbers, sort&compress]{natbib}
\usepackage[colorlinks=true,linkcolor=blue,citecolor=blue,backref=true,pagebackref=true]{hyperref}

\usepackage{esdiff}

\newtheorem{theorem}{Theorem}

\newtheorem{corollary}[theorem]{Corollary}

\newtheorem{remak}{Remark}

\newtheorem{exam}{Example}

\title{Directional derivatives of the singular values of matrices depending on several real parameters}
\author{Juan-Miguel Gracia\thanks{Department of Applied Mathematics and Statistics and O.R., University of the Basque Country UPV/EHU, Faculty of Pharmacy, Paseo de la Universidad, 7,
01006 Vitoria-Gasteiz, Spain,
\texttt{juanmiguel.gracia@ehu.es} }}
\date{May 14, 2020}

\begin{document}
\maketitle
\begin{abstract}
In this document I recapitulate some results by Hiriart-Urruty and Ye\cite{hiriart-urruty1995sensitivity} concerning the properties of differentiability and the existence of directional derivatives of the multiple eigenvalues of a complex Hermitian matrix function of several real variables, where the eigenvalues are supposed in a decreasing order. Another version of these results was obtained by Ji-guang Sun\citep{sun1988note,sun1988sensitivity,sun1990multiple}.
\end{abstract}

{
\hypersetup{linkcolor=blue}
\tableofcontents
}

\section{Differentiability   of the eigenvalues of a complex Hermitian matrix }

We will denote by $\s C$ the spectrum or set of eigenvalues of any complex square matrix $C$. Let $\Om$ be an open subset of $\R^p$ and let $A:\Om\to\Cnn$ be a matrix function of class $C^1$ such that for every $x\in\Om$ the matrix $A(x)$ is Hermitian, i.e. $A(x)^*=A(x)$ where $^*$ denotes the conjugate transpose. As it is well known the eigenvalues of $A(x)$ are real numbers; thus, there exist $n$ real functions defined on $\Om$, $\secsim \la n$, such that for all $x\in\Om$,
\[
\secudectres \la x n
\]
are the eigenvalues of $A(x)$. Let $\recoper m1n$; it is easy to prove that the function $\la_m\colon \Om \to \R$ is continuous. When the eigenvalue $\la_m(x_0)$ of $A(x_0)$ is simple, the function $\la_m$ is differentiable at $x_0\in\Om$. But in case of $\la_m(x_0)$ is a multiple eigenvalue of $A(x_0)$, $\la_m$ can be nondifferentiable at $x_0$. For example~\cite{kato:1982}, let
\[
A(x_1,x_2):=\begin{pmatrix}
 x_1 &  \mi x_2  \\
-  \mi x_2   & -x_1   
   \end{pmatrix}
\]
be for  $(x_1,x_2)\in\R^2$. It is obvious that for each $(x_1,x_2)\in\R^2$ the matrix $A(x_1,x_2)$ is Hermitian. Then 
\[
\begin{vmatrix}
\la-x_1 &  - \mi x_2 \\
  \mi x_2 & \la+x_1
\end{vmatrix}=\la^2-x_1^2-x_2^2;  
\] 
hence the eigenvalues of $A(x_1,x_2)$ are $\pm\sqrt{x_1^2+x_2^2}$. Observe that the matrix 
\[
A(0,0)=\begin{pmatrix}
 0 & 0   \\
 0  & 0
   \end{pmatrix}
\]
has a double eigenvalue; but neither the function $\la_1(x_1,x_2)=\sqrt{x_1^2+x_2^2}$,  nor
the function  $\la_2(x_1,x_2)=-\sqrt{x_1^2+x_2^2}$ are differentiable at $(0,0)$.

Let $d\in\R^p$ be a unitary vector, i.e. $\norma{d}_2=1$, where $\norma{\cdot}_2$ denotes the Euclidean norm. The directional derivative of the function $\la_m$ at the point $x_0$ with respect to $d$ is defined as the limit
\[
\la'_m(x_0,d):=\lim_{t\to 0^+}\frac{\la_m(x_0+td)-\la_m(x_0)}{t}
\]
whenever this limit exists.

Based on techniques and results from convex and nonsmooth analysis (in Clarke's sen\-se), Hiriart-Urruty and Ye proved  Theorems~\ref{th:10} and~\ref{th:20}. See~\cite[Theorem 4.5]{hiriart-urruty1995sensitivity}.

\begin{theorem}\label{th:10}
For all $x_0\in\Om$, for all unitary vector $d\in\R^p$, and for all $\recoper m1n$, there exists always
\[
\la'_m(x_0,d).
\]
\end{theorem} 

Moreover, it can be proved that $\la'_m(x_0,d)$ is equal to a determined eigenvalue
of a matrix constructed from $A(x_0)$ and $d$ in the following way: For each $x_0\in\Om$, there is a unitary matrix $U=[\secsim un]$ such that
\[
U^*A(x_0)U=\diag\big(\sectres \la{x_0}n\big).
\] 

Suppose that $\la_m(x_0)$ is a multiple eigenvalue of $A(x_0)$, of multiplicity $r_m$. Introduce two integers $i_m\ge 1,j_m\ge 0$ to precise the position that $\la_m(x_0)$ occupies  among the $r_m$ repeated eigenvalues that are equal to it.
Consider the detailed arrangement of the eigenvalues of $A(x_0)$:
\begin{multline*}
 \secudecuatro {\la}{1}{x_0}{m-i_m} > \SecuIgualCuatro  {\la}{m-i_m+1}{x_0}{m}       \\
=	\SecuIgualCuatro  {\la}{m+1}{x_0}{m+j_m}>	\secudecuatro {\la}{m+j_m+1}{x_0}{n}
\end{multline*}
That is to say, $j_m$ is the number of eigenvalues placed after the subscript $m$ that are equal to $\la_m(x_0)$; whereas $i_m$ is the number of eigenvalues placed before $m$ that are equal to $\la_m(x_0)$, plus one (we put  $\la_m(x_0)$ in this list). Hence, $j_m$ may be zero, $i_m\ge 1$, and $i_m+j_m=r_m$. When $m=1$, i.e. if we are considering $\la_1(x_0)$, we have $i_1=1,j_1=r_1-1$. When $m=n$, i.e. for $\la_n(x_0)$, we have $i_n=r_n, j_n=0$. In case $\la_m(x_0)$ is a simple eigenvalue,
$i_m=1,j_m=0$. Although the notation does not indicate it, the numbers $i_m,j_m$ and $r_m$ depend on $x_0$.

Let $U_2$ be the $n\x r_m$ matrix formed by the $(m-i_m+1)$th, \ldots, $(m+j_m)$th columns of the matrix $U$:
\[
U_2:=\left[u_{m-i_m+1},\ldots,u_{m+j_m}\right];
\]
i.e. $U_2$ is formed by $r_m$ orthonormal eigenvectors associated with the eigenvalue $\la_m(x_0)$ of $A(x_0)$. For each $\recoper j1p$ define 
\[
\frac{\partial A}{\partial x_j}(x_0)=\left(\frac{\partial a_{ik}}{\partial\, x_j}(x_0)\right)
\]
$a_{ik}(x)$ being the entries of $A(x)$.
We will call $F'(d)$ to the $r_m\x r_m$ matrix 
\[
F'(d):=U_2^* \left(\sum_{j=1}^p d_j\frac{\partial A}{\partial x_j}(x_0)\right)U_2
\]
for every unitary vector $d=(\secsim dp)\in\R^p$. Given that 
\[
 \overline{\left(\parcial {a_{ik}}{x_j}\right)}=\parcial {\bar{a}_{ik}}{x_j} =\parcial {a_{ki}}{x_j},
\]
we have that the matrix $\parcial A{x_j}$ is Hermitian, and so is $F'(d)$; indeed,
\begin{multline*}
 F'(d)^*=U_2^*\left(\sum_{j=1}^p d_j \parcial A{x_j}(x_0)\right)^*U_2            \\
		=U_2^*\left(\sum_{j=1}^p d_j \left[\parcial A{x_j}(x_0)\right]^*\right)U_2 	 \\
		=U_2^*\left(\sum_{j=1}^p d_j \parcial A{x_j}(x_0)\right)U_2 =F'(d).
\end{multline*}
Therefore, the eigenvalues of $F'(d)$ are real numbers. 

\begin{theorem}\label{th:20}
The directional derivative $\la'_m(x_0,d)$ is given by
\[
\la'_m(x_0,d)=\mu_{i_m}\big(F'(d)\big)
\]
where $\mu_{i_m}\big(F'(d)\big)$ is the $i_m$\emph{th} eigenvalue of $F'(d)$ when the eigenvalues are arranged in a decreasing order:
\[
\secudecuatro \mu {1}{F'(d)}{r_m}.
\]
\end{theorem}
A theorem related to Theorem~\ref{th:20} was proved by Ji-guang Sun~\cite[Theorem 3.1]{sun1990multiple} applying the implicit function theorem and the Rellich theorem.

\begin{theorem}\label{th:30}
The function
\[
t_m(x):=\la_{m-i_m+1}(x)+\cdots +\la_m(x)+\cdots +\la_{m+j_m}(x),\quad x\in\Om
\]
is differentiable at $x_0$.
\end{theorem}
See~\cite[Corollary 4.3]{hiriart-urruty1995sensitivity} for a proof of this theorem.

\begin{corollary}\label{cor:40}
There exists a neighborhood $V$ of $x_0$, $V\subset \Om $, in which the function 
\[
t_m(x):=\la_{m-i_m+1}(x)+\cdots +\la_m(x)+\cdots +\la_{m+j_m}(x)
\]
is differentiable.
\end{corollary}
\demo Let $V\subset \Om$ be a neighborhood of $x_0$,  sufficiently small  so that the inequalities
\[
\la_{m-i_m}(x)>\la_{m-i_m+1}(x), \quad \la_{m+j_m}(x)>\la_{m+j_m+1}(x)
\] 
hold when $x\in V$. 
Let $x_1$ be any point of $V$. Then the arrangement of the eigenvalues of $A(x_1)$
\[
\secudecuatro \la{m-i_m+1}{x_1}{m+j_m}
\]
may have groups of equalities. In view of Theorem~\ref{th:30}, the sum of the functions $\la_i$ corresponding to each one of these groups, is differentiable at $x_1$; therefore, as $t_m$ is the sum of these sums, we deduce that $t_m$ is differentiable at $x_1$. \hfill\qed

\section{Differentiability of the singular values of a complex matrix}

\textbf{Notation.} Let $m,n$ be positive integers and let $q:=\min(m,n)$. Given $B\in\Cmn$, let
\[
\secudectres \si{B} {q}
\]
be the singular values of $B$. For each $\recoper k1q$, it is said that a pair of vectors of \textbf{unit length} $y_k\in\Cmuno$, $ z_k \in \Cnuno$ are left and right
\textbf{singular vectors }of $B$ associated with the singular value $\si_k(B)$ if $Bz_k=\si_k(B)y_k$
and $B^*y_k=\si_k(B)z_k$.

\bigskip

 Let $A\colon \Om \to \Cmn$ be a matrix function of class $C^1$. For each $x\in \Om\subset \R^p$,
 let
 \[
 \secudecuatro s 1xq, 
\quad \text{  with } q:=\min(m,n),
 \]
be the singular values of the matrix $A(x)$ arranged in a decreasing order. Thus, we can define $q$ functions $s_i\colon \Om\to \R, \; \recoper i1q$. We are going to establish the properties of differentiability of these functions. By  Wielandt's lemma, the $m+n$ eigenvalues of the Hermitian matrix
\[
M(x):=\begin{pmatrix}
 0 & A(x)    \\
 A(x)^*  & 0  
   \end{pmatrix}\in \Cmmasn
\] 
are
\[
\secudecuatro s 1xq \ge 0= \cdots = 0 \ge \secudecuatro {-s}qx1
\]
(it may have repeated intermediate zeros), for all $x\in \Om$. Hence, the analogous results  to Theorems~\ref{th:10},~\ref{th:20} and~\ref{th:30}   for Hermitian matrices  are true.
\begin{theorem}\label{th:50}
Let $\recoper k1q$, $x_0\in\Om$, and $d\in \R^p$ be a unitary vector. Then there exists the directional derivative 
\[
s'_k(x_0,d).
\]
\end{theorem}

\bigskip

 Let $u\in \Cmuno, v\in \Cnuno$, where $u\neq 0$ or $v\neq 0$. Then
\[
\begin{pmatrix}
 u   \\
 v
   \end{pmatrix}
\]
is an eigenvector of 
\[
H:=\begin{pmatrix}
 0 & B   \\
 B^*  & 0 
   \end{pmatrix}
\]
 associated with its eigenvalue $\si_k(B)$  if and only if
\begin{equation}\label{eq:uno}  
Bv   = \si_k(B)u,   
\end{equation}   
\begin{equation} \label{eq:dos}  
 B^*u     =\si_k(B)v. 
\end{equation}
So, if $(y_k,z_k)\in \Cmuno\x \Cnuno$ is a pair of (left,right)-singular vectors of $B$ associated with the singular value $\si_k(B)$, then
\[
\begin{pmatrix}
y_k  \\ 
z_k 
   \end{pmatrix}
\]
is an eigenvector of $H$ corresponding to its eigenvalue $\si_k(B)$.
\bigskip

Let $x_0\in\Om$ be a fixed point, and let $W\in\Cmmasn$ a unitary matrix that diagonalizes $M(x_0)$:
\begin{equation}\label{eq:tres}
W^*M(x_0)W= \begin{pmatrix}
 s_1(x_0) &  &  & & & & 0 \\
    & \ddots & & & & &\\
 & & s_k(x_0) & &&&  \\
  &      & &\ddots & && \\
  & &        &  & -s_k(x_0) && \\
  & &           &  &  &\ddots &\\
 0 & &           &      &  & & -s_1(x_0)  
\end{pmatrix}
\end{equation}
Suppose that
\begin{multline*}
\secudecuatro s 1{x_0}{k-i_k} > \SecuIgualCuatro s{k-i_k+1}{x_0}k \\ = \SecuIgualCuatro s {k+1}{x_0}{k+j_k} > \secudecuatro s {k+j_k+1}{x_0}q \ge \cdots \ge -s_1(x_0)
\end{multline*}
are the eigenvalues of $M(x_0)$, where $s_k(x_0)$ is a multiple eigenvalue of multiplicity $r_k=i_k+j_k$, $i_k$ being the number of eigenvalues equal to $s_k(x_0)$ placed before the rank $k+1$, and $j_k$ is the number of eigenvalues equal to $s_k(x_0)$ situate after the rank $k$.

Call $W_2$ to the $(m+n)\x r_k$ matrix formed by the $(k-i_k+1)$th,\ldots,$(k+j_k)$th columns of the matrix $W$.  For each unitary vector $d=(\secsim dp)\in\R^p$, define
\[
F'(d):= W_2^* \left(\sum_{j=1}^p d_j \begin{bmatrix}
 O &  \parcial A{x_j}(x_0)  \\
 \left(\parcial A{x_j}(x_0)\right)^*  &   O
   \end{bmatrix} \right)W_2,
\]
which is an $r_k\x r_k$ Hermitian matrix. Then, by Theorem~\ref{th:20}, we have the next result.
\begin{theorem}\label{th:60}
For each unitary vector $d=$ $(\secsim dp)$$\in \R^p$
\[
s_k'(x_0,d)=\mu_{i_k}\big(F'(d) \big),
\]
$\mu_{i_k}\big(F'(d) \big)$ being the $i_k$\emph{th} eigenvalue of the matrix $F'(d)$
when we arrange the eigenvalues of this matrix in a decreasing order.
\end{theorem}

To facilitate the writing let $W_2$ be partitioned thus:
\[
W_2=\begin{bmatrix}
U_2\\
V_2
\end{bmatrix}
\]
where 
\[
U_2:=\left[\secuencia u{k-i_k+1}{k+j_k}\right],\quad V_2:=\left[\secuencia v{k-i_k+1}{k+j_k}\right],
\]
 \hfill $U_2\in \C^{m\x r_k}, V_2\in \C^{n\x r_k}$.
\begin{corollary}\label{cor:70}
For each unitary vector $d=$ $(\secsim dp)$$\in \R^p$
we have
\[
s_k'(x_0,d)=\mu_{i_k},
\]
where $\mu_{i_k}$ is the $i_k$\emph{th} eigenvalue of
\[
U^*_2 \left(\sum_{j=1}^p d_j\parcial A{x_j}(x_0)\right)V_2+\left(U^*_2 \left(\sum_{j=1}^p d_j\parcial A{x_j}(x_0)\right)V_2\right)^*,
\]
when the eigenvalues are ranked in a decreasing order.
\end{corollary}
\demo  Given that 
\[
W_2:=\begin{bmatrix}
U_2\\
V_2
\end{bmatrix},
\]
the matrix $F'(d)$ is given by
\begin{multline}
[U^*_2,V^*_2]\left(\sum_{j=1}^p d_j\begin{bmatrix}
 O & \parcial A{x_j} (x_0)   \\
 \left[\parcial {A}{x_j} (x_0)\right]^*  &   O
   \end{bmatrix}\right)\,\begin{bmatrix}
 U_2   \\
 V_2
   \end{bmatrix}\\
   =\sum_{j=1}^p d_j\left[V^*_2\left[\parcial{A}{x_j}(x_0)\right]^*,U^*_2\parcial{A}{x_j}
(x_0)\right]\;\begin{bmatrix}
 U_2   \\
 V_2
   \end{bmatrix}\\
   =V^*_2\sum_{j=1}^p d_j\left[\parcial{A}{x_j}(x_0)\right]^*\;U_2+U^*_2\sum_{j=1}^p d_j\parcial{A}{x_j}
(x_0)\;V_2\\
=U^*_2\sum_{j=1}^p d_j\parcial{A}{x_j}
(x_0)
\;V_2+\left(U^*_2\sum_{j=1}^p d_j\parcial{A}{x_j}
(x_0)\;V_2\right)^*.
\end{multline}
\hfill \cqd

\bigskip

 The sum of all singular values that coalesce with $s_k(x_0)$ at $x_0$
 is differentiable at $x_0$. Even more it is true as we can see in the next theorem.  
\begin{theorem}\label{th:80}
The function
\[
t_k(x):=\SecuSumaCuatro s{k-i_k+1}xk +\cdots +s_{k+j_k}(x)
\]
 is differentiable in a neighborhood $V\subset \Om$ of $x_0$.

\end{theorem}

The neighborhood $V$ is determined by the $x\in\Om$ sufficient close to $x_0$ in order that the inequalities 
\[
s_{k-i_k}(x)>s_{k-i_k+1}(x) \quad \text{  and  }\quad s_{k+j_k}(x)>s_{k+j_k+1}(x)
\]
hold.

\bigskip

From Corollary~\ref{cor:70} we can
give another description of $s'_k(x_0,d)$ in terms of singular vectors of $A(x_0)$ associated with $s_k(x_0)$. 
\begin{theorem}\label{th:30abril2010-10}
With the previous notation, let 
\[
Y=[y_{k-i_k+1},\ldots,y_{k+j_k}]\in \C^{m \x r_k},\quad Z=[z_{k-i_k+1},\ldots,z_{k+j_k}]\in\C^{n\x r_k}
\]
be matrices of orthonormal columns and such that $(y_\ell,z_\ell)\in \Cmuno\x \Cnuno$ is a pair of (left,right)-singular vectors of $A(x_0)$ associated with the singular value $s_k(x_0)$ for  $\recoper {\ell}{k-i_k+1}{k+j_k}$. 
Then $s'_k(x_0,d)$ is equal to the $i_k$\emph{th} eigenvalue of the $r_k\x r_k$ Hermitian matrix 
\[
G:=\dfrac{1}{2}\left[Y^*\left(\sum_{j=1}^p d_j\diffp{A}{{x_j}}(x_0)\right) Z+ Z^*\left(\sum_{j=1}^p d_j\left[\diffp{A}{{x_j}}(x_0)\right]^*\right) Y      \right],
\]
when the eigenvalues are ranked in a decreasing order.
\end{theorem}

\section{Function of Ikramov-Nazari} 
With the notations of~\cite{ikramov2003distance}, let
$(\xi_1,\xi_2,\xi_3,\xi_4)\in\R^4$, $A\in\Cnn$. Define
\[
Q(\xi_1,\xi_2,\xi_3,\xi_4):=\begin{pmatrix}
 A & \xi_1 I  & (\xi_3+\mi\,\xi_4)I  \\
 0  & A  & \xi_2 I\\
 0  &  0  & A
   \end{pmatrix},\quad n\ge 3.
\]
Set 
\[
f(\xi):=s_{3n-2}\big(Q(\xi)\big).
\]
Suppose that the function $f$ attains a local maximum at a given $\xi^0\in\R^4$, say
$s_0:=s_{3n-2}\big(Q(\xi^0)\big)$. Let us also assume that $s_0>0$ and it is a multiple singular
of $Q(\xi^0)$. With the above notations, there are
$i_{3n-2}$ singular values before the place $3n-2+1$ and $j_{3n-2}$ singular values  after the place $3n-2$ equal to $s_{3n-2}\big(Q(\xi^0)\big)$. 
To shorten notation, we let $p$ and $q$ stand for $i_{3n-2}$ and $j_{3n-2}$, respectively.  
Thus, the multiplicity of $s_0$ is
$m=p+q$. Hence,
\begin{multline*}
\secudecuatro{s}1{Q(\xi^0)}{3n-2-p}\\ > \SecuIgualCuatro{s}{3n-2-p+1}{Q(\xi^0)}{3n-2}\\
= \SecuIgualCuatro{s}{3n-2+1}{Q(\xi^0)}{3n-2+q}\\ >           
	\secudecuatro{s}{3n-2+q+1}{Q(\xi^0)}{3n}.		 
\end{multline*}
Here $p\ge 1$ and $q\ge 0$. The function
\[
t(\xi):=\SecuSumaCuatro{s}{3n-2-p+1}{Q(\xi)}{3n-2+q}
\]
is differentiable in a neighborhood of $\xi^0$. Also for each $\recoper k1{3n}$ and each unitary vector $d\in\R^4$, the function 
\[
g_k(\xi):=s_k\big(Q(\xi)\big)
\]
admits the directional derivative
\[
g'_k(\xi^0,d).
\]
Observe that the used notation implies
\[
f(\xi)=g_{3n-2}(\xi),\quad \xi\in\R^4.
\]
Next, we determine the relationship between the directional derivatives $f'(\xi^0,d)$ and $f'(\xi^0,-d)$.
Given that $f$ has a local maximum at $\xi^0$, it follows that for all $e\in\R^4$,
\[
f'(\xi^0,e):=\lim_{h\to 0^+}\frac{f(\xi^0+he)-f(\xi^0)}{h}\le 0.
\]
Thus, $f'(\xi^0,d)\le 0$ and $f'(\xi^0,-d)\le 0$. 
What conditions must be satisfied in order for $f'(\xi^0,d)=0$ to hold for all unit vector $d\in\R^4$?
By Theorem~\ref{th:60},
 $f'(\xi^0,d)$ is equal to  $\mu_p(d)$, $p$th eigenvalue of the $m\x m$ matrix
\[
F'(d)=[U^*_2,V^*_2]\left(\sum_{j=1}^4 d_j\begin{bmatrix}
 0 & \parcial Q{\xi_j} (\xi^0)   \\
 \parcial {Q^*}{\xi_j} (\xi^0)  &   0
   \end{bmatrix}\right)\begin{bmatrix}
 U_2   \\
 V_2
   \end{bmatrix}
\]
with
\begin{align*}
U_2    & =  [\secuencia u{3n-2-p+1}{3n-2+q}]         \\
V_2    & =  [\secuencia v{3n-2-p+1}{3n-2+q}]         
\end{align*}
where $u_j$ and $v_j$ are the left and right singular vectors
\[
\left. 
\begin{matrix}
Q(\xi^0)v_j=s_0u_j\\
Q(\xi^0)^*u_j=s_0v_j
\end{matrix}\right\}\qquad \recoigu j{3n-2-p+1}{3n-2+q},
\] 
and the eigenvalues of $F'(d)$ are arranged in this way
\begin{equation}\label{eq:70} 
\secudecuatro \mu 1dp \ge \cdots \ge \mu_m(d)
\end{equation}
Therefore, 
\[
f'(\xi^0,d)=\mu_p(d).
\]
Again by Theorem~\ref{th:60}, we deduce that $f'(\xi^0,-d)$ is equal to the $p$th eigenvalue of the Hermitian matrix
$F'(-d)$. 
But, it is worth noting that $f'(\xi^0,-d)$ is not necessarily equal to $\mu_p(-d)$. In fact, if
\[
\secdec \al m
\]
are the eigenvalues of $F'(-d)$, then
\[
f'(\xi^0,-d)=\al_p.
\]
As $F'(-d)=-F'(d)$, it follows 
\begin{equation}\label{eq:80} 
\secudecuatro {-\mu}mdp \ge \cdots \ge -\mu_1(d)
\end{equation}
are the eigenvalues of $F'(-d)$; whence,
\begin{equation}\label{eq:90} 
f'(\xi^0,-d)=\al_p=-\mu_{m-(p-1)}(d).
\end{equation}
Now it is necessary to analyze the relative positions of the indices $p$ and $m-(p-1)$. 

If $p\le m-(p-1)$, then $\mu_p(d)\le 0$, and it follows that
\[
0\ge \secudecuatro \mu pd{m-(p-1)}\ge \cdots \ge \mu_m(d).
\]
Hence, $0\ge \mu_{m-(p-1)}(d)$ and therefore $\al_p = -\mu_{m-(p-1)}(d)\ge 0$, but $\al_p=f'(\xi^0,-d)\le 0$. Thus, $\al_p=0$; i.e. $f'(\xi^0,-d)=0$. Given that $f$ has a local maximum at $\xi^0$, for all unit vector $e\in\R^4$, we have
\[
f'(\xi^0,e)=0.
\]
 Doutbful   case: If $p>m-(p-1)$,  then $\mu_{m-(p-1)}(d)\ge \mu_p(d)$. But, although  $\mu_p(d)\le 0$, it is not guaranteed that the inequality $\mu_{m-(p-1)}(d)\le 0$ holds.
 
\section{Average of singular values}
We know that the average of singular values of $Q(\xi)$ that coalesce with the $m$-multiple singular value
$s_{3n-2}\big(Q(\xi^0)\big)$ at $\xi=\xi^0$, is a differentiable function
in a neighborhood of $\xi^0$. Thus we consider the differentiable function
\[
H(\xi):=t(\xi)-ms_0;
\]
obviously, $H(\xi^0)=0$. Hence, the point $\xi^0$ belongs to the level hypersurface
of level $0$ of the function $H(\xi)$. Let
\[
\nabla H(\xi^0)=\left(\parcial{H}{{\xi}_1}(\xi^0),\parcial{H}{{\xi}_2}(\xi^0),
\parcial{H}{{\xi}_3}(\xi^0),\parcial{H}{{\xi}_4}(\xi^0)\right)
\]
be the gradient of $H(\xi)$ at $\xi^0$. Let $d\in\R^4$ such that
\[
\nabla H(\xi^0)\cdot d=0,
\] 
where $\cdot$ denotes the ordinary scalar product in $\R^4$. Then, by the chain rule,
\[
H'(\xi^0,d)=\nabla H(\xi^0)\cdot d=0.
\]
This implies
\[
0=g'_{3n-2-p+1}(\xi^0,d)+\cdots +g'_{3n-2}(\xi^0,d)+\cdots +g'_{3n-2+q}(\xi^0,d);
\]
if we consider the $m\x m$ Hermitian matrix $F'(d)$, it means that the sum of its eigenvalues is zero:
\[
0=\mu_1(d)+\cdots +\mu_p(d)+\cdots +\mu_m(d).
\]

When $p=1$, this is equivalent to say that $s_{3n-2}\big(Q(\xi^0)\big)$ is the first value of the chain of singular values equal to $s_0$, then all the functions
\[
g_{3n-2}(\xi),g_{3n-2+1}(\xi),\ldots,g_{3n-2+q}(\xi)
\]
 take the same value  at $\xi^0$, and it is equal to $s_0$. Moreover, all these functions have at $\xi^0$ a local maximum, because of
\[
f(\xi):=g_{3n-2}(\xi)\ge g_{3n-2+1}(\xi)\ge \cdots \ge g_{3n-2+q}(\xi).
\]
This implies that for all unitary $d\in\R^4$,
\[
\forall \recoigu k{3n-2}{3n-2+q},\quad g'_k(\xi^0,d)\le 0;
\]
therefore, $\mu_1(d)\le 0,\ldots,\mu_m(d)\le 0$, and, given that $t(\xi)$ has a local maximum at $\xi^0$ and is differentiable at $\xi^0$, we have
\[
\nabla t(\xi^0)=\n{0};
\]
whence $\nabla H(\xi^0)=\n{0}$ and for all $\recoigu k{3n-2}{3n-2+q}$, $g_k'(\xi^0,d)=0$; in particular, $f'(\xi_0,d)=0$. This is proved because
$0=\SecuSumaCuatro \mu 1 d m$; since $\forall k,\mu_k(d)\le 0$, we obtain $\forall k, \mu_k(d)=0$; consequently, $\forall k, g'_k(\xi^0,d)=0$.

From now on let $p$ be any integer from the range we are considering. Furthermore, suppose that for all $k=3n-2-p+1,\ldots,3n-2+q$, all the functions $g_k(\xi)$ have a local maximum at $\xi^0$. Then for all unitary $d
\in \R^4$,$g'_k(\xi^0,d)\le 0$. As $t(\xi)$ has a local maximum at $\xi^0$, $t'(\xi^0,d)=0$; but
\[
t'(\xi^0,d)=g'_{3n-2-p+1}(\xi^0,d)+\cdots+g'_{3n-2+q}(\xi^0,d);
\]
consequently, $f'(\xi^0,d)=0$.

When some of the functions $g_k(\xi)$ have a local maximum at $\xi^0$ and any others have a local minimum at $\xi^0$, the analysis becomes more complicated and I do not obtain any conclusion.

\section{Remark}
In January 31, 2005, I wrote an e-mail to J.B. Hiriart-Urruty asking him whether his results in \cite{hiriart-urruty1995sensitivity} for real symmetric matrices could be generalized to complex Hermitian matrices. He forwarded my message to M. Torki \cite{torki:2005}, who answered  affirmatively. Moreover, Torki  told me that his results in \cite{torki:2001} for second order directional derivatives and real symmetric matrices were also true for the Hermitian case. 

A particular case of  Theorem~\ref{th:30abril2010-10} about the right derivative of the function $t\mapsto \si_k(A+tB)$ at $t=0$, where $t$ is real and $A,B$ are $n\x n$ complex matrices,  was obtained by Lippert~\cite[Lemma A.5]{lippert2005fixing} by a different  method. See also Corollary 11 in~\cite{armentia2018perforated}.

\subsection*{Acknowledgments}

I thank Gorka Armentia and Francisco E. Velasco for the talks we had about this  topic.

\end{document}